\documentclass[11pt]{amsart}
\usepackage{hyperref,flafter,ifsym,wasysym,color}
\usepackage{epsfig,bbm}
\usepackage[english]{babel}
\usepackage{amsfonts}
\usepackage{amssymb}
\usepackage{amsmath}
\usepackage{amsthm}
\usepackage{stmaryrd}
\usepackage{latexsym}
\usepackage{amscd}
\usepackage{epsf}
\usepackage{caption}
\usepackage[lmargin=1in,rmargin=1in,tmargin=0.8in,bmargin=0.8in]
{geometry}
\usepackage{enumitem}
\usepackage{todonotes}
\usepackage{tikz,tikz-cd}
\usetikzlibrary{positioning}
\usetikzlibrary{patterns}
\usetikzlibrary{shapes}
\usetikzlibrary{plotmarks}
\newcommand{\X}{\mathsf{X}}

\newcommand{\Os}{\mathsf{O}}

\newcommand{\ra}{\rightarrow}

\newcommand{\Mod}{\Mcal}
\newcommand{\Mhat}{\widehat{\Mod}}
\newcommand{\Nhat}{\widehat{\Ncal}}
\newcommand{\sig}{\sigma}
\newcommand{\Group}{\mathsf{G}}

\newcommand{\states}{\mathsf{S}}

\newcommand{\emptyrect}{\mathrm{Rect}^\circ}
\newcommand{\CFT}{\mathrm{CF}}

\newcommand{\Dbb}{\mathbb{D}}

\newcommand{\Obb}{\mathbb{O}}

\newcommand{\Xbb}{\mathbb{X}}

\newcommand{\Hcal}{\mathcal{H}}

\newcommand{\Mcal}{\mathcal{M}}
\newcommand{\Ncal}{\mathcal{N}}

\newcommand{\Scal}{\mathcal{S}}

\newcommand{\ofrak}{\mathfrak{o}}

\newcommand{\ufrak}{\mathfrak{u}}
\newcommand{\vfrak}{\mathfrak{v}}

\newtheorem{thm}{Theorem}[section]

\newtheorem{remark}[thm]{Remark}

\theoremstyle{definition}

\def\endproof{\relax\ifmmode\expandafter\endproofmath\else
  \unskip\nobreak\hfil\penalty50\hskip.75em\hbox{}\nobreak\hfil\bull
  {\parfillskip=0pt \finalhyphendemerits=0 \bigbreak}\fi}
\def\endproofmath$${\eqno\bull$$\bigbreak}
\def\bull{\vbox{\hrule\hbox{\vrule\kern3pt\vbox{\kern6pt}\kern3pt
\vrule}\hrule}}

\newcommand{\R}{\mathbb{R}}

\newcommand{\Z}{\mathbb{Z}}

\newcommand\x{\mathbf x}
\newcommand\y{\mathbf y}

\newcommand\alphas{\mbox{\boldmath$\alpha$}}

\newcommand\betas{\mbox{\boldmath$\beta$}}

\newcommand\SpinC{\mathrm{Spin}^c}

\begin{document}
\title{On sign assignments in link Floer homology}%
\author{Eaman Eftekhary}%
\address{School of Mathematics, Institute for Research in 
Fundamental Sciences (IPM), P. O. Box 19395-5746, Tehran, Iran}%
\email{eaman@ipm.ir}
\begin{abstract}
In this short note, we compare the combinatorial sign assignment 
of Manolescu, Ozsv\'ath, Szab\'o and Thurston for grid homology 
of knots and links in $S^3$ with the sign assignment coming from 
a coherent system of orientations on  Whitney disks. 
Although these constructions produce 
different signs, a small modification of the convention in either
of the two methods results in identical sign assignments, and thus 
identical chain complexes.  
\end{abstract}
\maketitle
\section{Introduction}
Extending the construction of Heegaard Floer invariants for 
three-manifolds \cite{OS-3m1}, knot Floer homology was introduced by 
Ozv\'ath and Szab\'o \cite{OS-knot}, and independently by Rasmussen
\cite{Ras1}, c.f. \cite{Ef-LFH} for the case of homologically non-trivial 
knots. The construction was refined in \cite{OS-link} to construct 
invariants of pointed links, which were defined over $\Z/2\Z$.
 Manolescu, Ozsv\'ath and Sarkar gave 
a combinatorial description of these link invariants \cite{MOS},
which was further explored by Manolescu, Ozsv\'ath, Szab\'o 
and Thurston \cite{MOST} to present the combinatorially defined 
{\emph{grid homology}} for knots and links in $S^3$ with integer
coefficients. On the other hand, Alishahi and the author \cite{AE-1}
refined sutured Floer homology of Juh\'asz \cite{Juh} (hat version 
of Heegaard Floer homology for sutured manifolds) to a theory that is 
defined over integers and includes other versions of Heegaard 
Floer homology for $3$-manifolds and pointed links as special cases.
The naturality of the construction is addressed by Alishahi and 
the author \cite{AE-2}, and independently by Zemke \cite{Ian}, but 
only over $\Z/2\Z$.\\

The goal of this short note is to compare the sign assignment of  
\cite{MOST} with the sign assignment of \cite{AE-1}. The signs in 
the construction of \cite{MOST} are defined combinatorially, 
by assigning a value in 
$\{-1,+1\}$ to every {\emph{empty rectangle}} which contributes 
to the differential of the complex. On the other hand,
the construction of \cite{AE-1} describes how the 
moduli spaces of holomorphic curves appearing in link Floer homology
 may be equipped with a
{\emph{coherent system of orientations}}, denoted by $\ofrak$. 
The key is the correct choice of orientation on the moduli spaces 
of boundary degenerations. In \cite{AE-1}, $\ofrak$
is chosen so that the signed count of $J$-holomorphic representatives 
of every boundary degeneration class with Maslov index $2$ is 
$+1$.  The signs assigned to holomorphic disks by such $\ofrak$, 
which are called  coherent systems of orientations with {\emph{ 
positive boundary degenerations}} is {\emph{not}} compatible with 
the sign assignment of \cite{MOST}. Alternatively, one may  
consider coherent systems of orientations with {\emph{negative 
boundary degenerations}, i.e. those $\ofrak$ such that the signed 
count of $J$-holomorphic representatives of every boundary 
degeneration class with Maslov index $2$, is $-1$. Then the two 
conventions for sign assignment agree: Every grid diagram is also 
a Heegaard diagram and every empty rectangle corresponds to the 
homotopy class of a Whitney disk. Since the moduli space 
of holomorphic representatives of every empty rectangle consists 
of a single point, a coherent system of orientations $\ofrak$ 
gives a map $\Scal_\ofrak:\emptyrect\ra \{-1,+1\}$,
where $\emptyrect$ denoted the set of all empty rectangles.
\begin{thm}\label{prop:main}
For every coherent system of orientations $\ofrak$ with 
negative boundary degenerations on a 
grid diagram $\Dbb$ for a link $L\subset S^3$, the function
$\Scal_\ofrak:\emptyrect\ra \{-1,+1\}$
is a true sign assignment in the sense of 
\cite[Definition 4.1]{MOST}. Moreover, every true sign assignment 
is of the form $\Scal_\ofrak$ for a uniquely determined 
coherent system of orientations $\ofrak$ with negative boundary 
degenerations. In particular, the grid chain 
complex associated with $\Dbb$ and the true sign assignment 
$\Scal_\ofrak$ is naturally identified with the link Floer complex
associated with $(\Dbb,\ofrak)$.  
\end{thm}

\section{Proof of the theorem}
Let us assume that $\Dbb$ is a grid diagram, representing some link 
$L\subset S^3$, which does not play a significant role for the 
purposes of this paper. Then $\Dbb$ consists of a torus $T$, 
a collection of $n$ horizontal circles 
$\alphas=\{\alpha_1,\ldots,\alpha_n\}$,
a collection of $n$ vertical circles 
$\betas=\{\beta_1,\ldots,\beta_n\}$, a collection of $\X$ markings
$\Xbb=\{\X_1,\ldots,\X_n\}$ and a collection 
$\Obb=\{\Os_1,\ldots,\Os_n\}$ of $\Os$ markings. We assume that 
horizontal and vertical circles appear with the (cyclic) order 
determined by their indices. Let us assume that 
$L$ has $m$ link components $L_1,\ldots, L_m$ and that $\X_i$ is on 
$L_{\imath(i)}$. Correspondingly, we obtain a Heegaard diagram 
\[\Hcal=(T,\alphas,\betas,\Xbb,\Obb).\]
The chain complex $\CFT(\Dbb)$, both in holomorphic curve approach 
and the combinatorial approach, is  freely generated over 
$\Z[\ufrak_1,\ldots,\ufrak_n,\vfrak_1,\ldots,\vfrak_m]$ by the 
{\emph{grid states}}. The grid states are in one to one 
correspondence with permutations $\sigma\in S_n$: associated with 
every $\sigma$ as above we have a generator 
$\x_\sigma=\{x_{1,\sigma(1)},\ldots,x_{n,\sigma(n)}\}$ where 
$x_{i,j}$ is the unique intersection of $\alpha_i$ and $\beta_j$. 
 The set of grid states is denoted by 
$\states=\states(\Dbb)$. For every $\x,\y\in\states$, 
$\emptyrect(\x,\y)$ denotes the set of empty rectangles which 
connect $\x$ to $\y$. It is clear that there is an inclusion  
$\emptyrect(\x,\y)\subset\pi_2(\x,\y)$. Associated with every 
$r\in\emptyrect(\x,\y)$, the moduli space $\Mhat(r)=\Mod(r)/\R$
of $J$-holomorphic representatives of $r$ (divided by the translation 
action of $\R$) consists of a single element. Implicitly, we are 
of course fixing a path $J$ of almost complex structures on $T$
throughout our discussion.\\

Let us quickly review the sign assignment coming from the 
construction of \cite{AE-1}. Since all grid states represent the 
unique $\SpinC$ structure on $S^3$, for every $\x,\y\in\states$, 
$\pi_2(\x,\y)$ is non-empty. For every $r\in\pi_2(\x,\y)$, the 
determinant line bundle associated with the linearization of 
time-dependent Cauchy-Riemann operator over the space of 
representative of  $r$ is trivial. One may choose a 
coherent system of orientations as follows. For $\sig\neq Id_{S_n}$,
fix the class $r_\sig$ of a Whitney disk in $\pi_2(\x_0,\x_\sig)$,
where $\x_0=\x_{Id_{S_n}}$. Moreover, let 
$a_1,\ldots,a_n\in\pi_2^\alpha(\x_0)$ denote the classes of $\alpha$
boundary degenerations, where $a_i$ corresponds to the thin 
cylinder bounded between $\alpha_i$ and $\alpha_{i+1}$ (with 
$\alpha_{n+1}=\alpha_1$). Similarly, let 
$b_1,\ldots,b_n\in\pi_2^\beta(\x_0)$ denote the classes of $\beta$
boundary degenerations, where $b_i$ corresponds to the thin 
cylinder bounded between $\beta_i$ and $\beta_{i+1}$.
The moduli spaces $\Ncal(a_i)$ and $\Ncal(b_j)$ of $J$-holomorphic 
representatives of $a_i$ and $b_j$ (respectively) 
are $2$-dimensional, and 
\[\Group=\left\{\left(\begin{array}{cc}a&b\\0&\frac{1}{a}\end{array}
\right)\ \Big|\ a\in\R^+, b\in\R\right\}<\mathrm{PSL}_2(\R)\]
acts on them. It was observed in \cite[Section 4]{AE-1} that the 
orientation on the determinant line 
bundle may be chosen so that the signed count of points in either
of $\Nhat(a_i)=\Ncal(a_i)/\Group$ and $\Nhat(b_j)=\Ncal(b_j)/\Group$ 
is $+1$. If the aforementioned choice of orientation is fixed 
over $a_i$ and $b_j$, the coherent system of orientations 
is said to have {\emph{positive boundary degenerations}}. 
Alternatively, as we will assume from here on, one may choose the 
orientation on the determinant line bundle so that the signed count 
of points in either of $\Nhat(a_i)$ and $\Nhat(b_j)$ is $-1$, and 
the coherent system of orientations is then said to have negative 
boundary degenerations. Note that each $a_i\in\pi_2^\alpha(\x_0)$ 
and $b_j\in\pi_2^\beta(\x_0)$ may also be considered as a Whitney 
disk in $\pi_2(\x_0,\x_0)$. In both cases,
the above choices of orientation is compatible with the equality
\[a_1\star a_2\star \cdots \star a_n=b_1\star b_2\star 
\cdots \star b_n\in\pi_2(\x_0,\x_0).\]

Either of the two possible orientations on the determinant line 
bundle over $r_\sig$ may be chosen by the coherent system $\ofrak$ of 
orientations (with negative boundary degenerations). Having fixed the 
above choices, $\ofrak$ picks a well-defined orientations on all 
Whitney disks as follows. Given  $r\in\pi_2(\x_\sig,\x_\tau)$, 
we have $r_\sig\star r=p\star r_\tau$, 
where $p\in\pi_2(\x_0,\x_0)$ is a juxtaposition of the $\alpha$ and 
$\beta$ boundary degenerations (viewed as classes of Whitney disks 
in $\pi_2(\x_0,\x_0)$). The orientation on the determinant line 
bundle over $p$ is determined  by the orientation on the determinant 
line bundles over $a_i$ and $b_j$. Thus, the orientation on the 
determinant line bundle over $r$, which is determined by the choices 
of orientation over $p$, $r_\sig$ and $r_\tau$, is uniquely 
determined by our earlier choices. Associated with each $a_i$, one 
finds an $\alpha$ boundary degeneration class 
$a_i^\sig\in\pi_2^\alpha(\x_\sig)$ which is defined by 
$a_i\star r_\sig=r_\sig\star a_i^\sig$. Either of the two choices 
of orientation on the determinant line bundle on $r_\sig$ induces 
the {\emph{negative}} orientation on $a_i^\sig$, i.e. the orientation 
with the property that the number of points in $\Nhat(a_i^\sig)$,
counted with sign, is $-1$. Similarly, one can define $b_i^\sig$,
and observe that the number of points in $\Nhat(b_i^\sig)$,
counted with sign, is $-1$.\\

Since for every $r\in\emptyrect$, $\Mhat(r)$ consists of a single 
point, the coherent system of orientations $\ofrak$  may be used 
to define a map \[\Scal_\ofrak:\emptyrect\ra \{-1,+1\}.\]
If $\ofrak'$ is another coherent system of orientations,
a function $f:S_n\ra \{-1,1\}$ describes the difference between 
$\ofrak$ and $\ofrak'$, in the sense that 
\[\Scal_{\ofrak'}(r)=f(\sig)\cdot f(\tau)\cdot \Scal_\ofrak(r)
\quad\quad\forall\ \ r\in\pi_2(\x_\sig,\x_\tau).\]
By \cite[Theorem 4.2]{MOST}, in order to prove Theorem~\ref{prop:main}
it suffices to prove the following statements:

\begin{itemize}
\item[(Sq)] For any four distinct $r_1,r_2,r_1',r_2'\in\emptyrect$ 
with $r_1\star r_2=r_1'\star r_2'$ we have
\[\Scal_\ofrak(r_1)\cdot\Scal_\ofrak(r_2)=-\Scal_\ofrak(r_1')
\cdot\Scal_\ofrak(r_2').\]
\item[(V)] If $r_1\in\emptyrect(\x_\sig,\x_\tau)$ and 
$r_2\in\emptyrect(\x_\tau,\x_\sig)$ satisfy
$r_1\star r_2=b_i^\sig$, then 
$\Scal_\ofrak(r_1)\cdot\Scal_\ofrak(r_2)=-1$.
\item[(H)] If $r_1\in\emptyrect(\x_\sig,\x_\tau)$ and 
$r_2\in\emptyrect(\x_\tau,\x_\sig)$ satisfy
$r_1\star r_2=a_i^\sig$, then 
$\Scal_\ofrak(r_1)\cdot\Scal_\ofrak(r_2)=+1$.\\
\end{itemize}

In order to prove (Sq), note that $r=r_1\star r_2$ has Maslov index 
$2$, and for generic $J$, $\Mhat(r)$ is a $1$-dimensional
manifold. The coherent system of orientations equips 
$\Mhat(r)$ with an orientation. The boundary of $\Mhat(r)$ 
correspond to possible degenerations of $r$. Since 
$r=r_1\star r_2=r_1'\star r_2'$ is decomposed in two ways, 
one concludes that no boundary degenerations appear as weak 
limits of $J$-holomorphic curves in $\Mhat(r)$.
Moreover, in every degeneration of $r$ to two Whitney disks, 
either of the Whitney disks is an empty rectangle. Nevertheless,
there are at most two such degenerations of $r$. In particular, 
the boundary of $\Mhat(r)$ is precisely
\[\left(\Mhat(r_1)\times\Mhat(r_2)\right)\coprod
\left(\Mhat(r_1')\times\Mhat(r_2')\right).\]
Since the determinant line bundle over $r$ is oriented compatible 
with the decompositions $r=r_1\star r_2$ and $r=r_1'\star r_2'$, it 
follows that 
\[\Scal_\ofrak(r_1)\cdot\Scal_\ofrak(r_2)+
\Scal_\ofrak(r_1')\cdot\Scal_\ofrak(r_2')=0.\]

To prove (V), consider the oriented $1$-manifold $\Mhat(r)$, where
$r=r_1\star r_2\in\pi_2(\x_\sig,\x_\sig)$. The boundary points
of $\Mhat(r)$ are either in correspondence with degenerations 
of $r$ to two Whitney disks of Maslov index $1$, or boundary 
degenerations.  Since $r$ may only be decomposed as $r_1\star r_2$,
the former type of boundary degenerations is in correspondence with 
$\Mhat(r_1)\times\Mhat(r_2)$. On the other hand, $r$ also corresponds
to a $\beta$ boundary degeneration $b_j^\sig\in\pi_2^\beta(x_\sig)$.
Correspondingly, we obtain boundary points which are in 
correspondence with $\Nhat(b_j^\sig)$. By \cite[Lemma 5.4]{AE-1},
the orientation of $\Nhat(b_j^\sig)$ induced by the coherent 
system of orientations $\ofrak$ is the opposite of the orientation 
induced on it as the boundary of $\Mhat(r)$. Since the signed count
of points in $\Nhat(b_j^\sig)$ with the orientation induced by 
$\ofrak$ is $-1$, it follows that 
\[\Scal_\ofrak(r_1)\cdot\Scal_\ofrak(r_2)-(-1)=0.\]
The proof of (H) is similar. The only difference is that 
\cite[Lemma 5.4]{AE-1} now implies that the orientation of 
$\Nhat(a_i^\sig)$ induced by $\ofrak$ agrees with the orientation
induced on it as the boundary of $\Mhat(r_1\star r_2)$, which gives
\[\Scal_\ofrak(r_1)\cdot\Scal_\ofrak(r_2)+(-1)=0.\]

\begin{remark}
Instead of using coherent systems of orientations with negative 
boundary degenerations, one can modify the convention on the 
combinatorial side by modifying (V) and (H), so that 
$\Scal_\ofrak(r_1)\cdot\Scal_\ofrak(r_2)=1$ when 
$r_1\star r_2=b_j^\sig$ and 
$\Scal_\ofrak(r_1)\cdot\Scal_\ofrak(r_2)=-1$ when 
$r_1\star r_2=a_i^\sig$. These {\emph{false}} sign assignments are 
in one-to-one correspondence with true sign assignments: if 
$\Scal_\ofrak$ is a true sign assignment, define 
\[\Scal'_\ofrak(r)=\mathrm{sgn}(\sig)\cdot \Scal_\ofrak(r)\quad
\quad\forall\ \ r\in\emptyrect(\x_\sig,\x_\tau).\]
Then $\Scal'_\ofrak$ is a false sign assignment. False sign 
assignments may be used in grid homology exactly like true sign 
assignments, to produce homology groups defined over integers, which
are knot/link invariants. This modification has the advantage that 
boundary degenerations are counted with positive sign, which 
is perhaps more natural.
\end{remark}
\bibliographystyle{hamsalpha}
\providecommand{\bysame}{\leavevmode\hbox to3em{\hrulefill}\thinspace}
\providecommand{\href}[2]{#2}
\providecommand{\eprint}{\begingroup \urlstyle{rm}\Url}

\end{document}